\documentclass[12pt,a4paper, twoside]{article}

\usepackage{amsmath,amssymb,amsthm,amsfonts,mathrsfs,amscd,environ}
\usepackage{dsfont}
\usepackage{latexsym,enumerate,color,geometry,extarrows}
\usepackage{xcolor}
\usepackage{verbatim,fancyhdr,enumitem}
 \NewEnviron{ews}{%
\begin{equation}\begin{split}
  \BODY
\end{split}\end{equation}
}

\NewEnviron{ews*}{%
\begin{equation*}\begin{split}
  \BODY
\end{split}\end{equation*}
}

 \newtheorem{thm}{Theorem}[section]
 \newtheorem{lem}[thm]{Lemma}
 
 \newtheorem{exa}[thm]{Example}
 
 \newtheorem{defn}{Definition}[section]
 \newtheorem{rem}{Remark}[section]
 
 \numberwithin{equation}{section}

\def\dif{{\mathord{{\rm d}}}}
\def\no{\nonumber}

\def\mR{{\mathbb R}}
\def\mE{{\mathbb E}}

\def\mN{{\mathbb N}}

\def\sF{{\mathscr F}}

\def\bd{\begin{defn}}
\def\ed{\end{defn}}
\def\bl{\begin{lem}}
\def\el{\end{lem}}
\def\bt{\begin{thm}}
\def\et{\end{thm}}
\def\br{\begin{rem}}
\def\er{\end{rem}}

\allowdisplaybreaks

\title{{\bf Stability of  Numerical Solution to    Pantograph Stochastic Functional Differential  Equations}\footnote{Supported in
 part by  NNSFC (61876192, 11626236) and  the Fundamental Research Funds for the Central Universities of South-Central University for Nationalities (CZY15017, KTZ20051, CZT20020).} }

\author{
{\bf Hao Wu$^{a)}$, Junhao Hu$^{a)}$, Chenggui Yuan$^{b)}$}\\
\footnotesize{$^{a)}$School of Mathematics and Statistics,
South-Central University for Nationalities}\\
\footnotesize{ Wuhan, Hubei 430000, P.R.China}\\
\footnotesize{Email: wuhaomoonsky@163.com},
\footnotesize{ junhaohu74@163.com}\\
\footnotesize{$^{b)}$Mathematic department, Swansea University, Bay campus, SA1 8EN, UK}\\
\footnotesize{Email: C.Yuan@Swansea.ac.uk}\\
}

\begin{document}

\maketitle

\begin{abstract}
In this paper, we  study the convergence of  the Euler-Maruyama numerical solutions for pantograph stochastic functional differential  equations which was proposed in \cite{why}.
 We also show  that the  numerical solutions  have the properties of  almost surely polynomial stability and  exponential stability  with the help of semi-martingale convergence theorem.
\end{abstract}\noindent

\noindent


\section{Introduction}

In order to solve a problem on the pantograph of an electric
locomotive, Ockendon and Tayler \cite{OT} proposed pantograph differential equations (PDEs). PDEs then have used in many areas such as electric dynamics,  dielectric materials and continuum mechanics. Additionally, many researchers have extended PDEs to   pantograph stochastic differential equations (PSDEs), so as to capture the practice problems more precisely. The form of PSDEs is as follows:
\begin{align*}
\begin{cases}
 \dif x(t)= b( x(t), x(q t), t)\dif t +\sigma( x(t),x(q t),t)\dif  B(t),\\
 x(0)=x_{0}.
 \end{cases}
 \end{align*}
where $q$  is a fixed constant satisfying $0<q<1.$  Many properties of PSDEs have been studied.  For example,      Milo\v{s}evi\'{c} \cite{LR10}  studied existence, uniqueness and  almost sure polynomial stability of solution to a class of highly nonlinear PSDEs and the Euler-Maruyama (EM) approximation. Shen et al. \cite{SFM} investigated the exponential stability of highly nonlinear neutral PSDEs. Guo et al. \cite{glhj} discussed the stability of numerical solutions for the PSDEs with variable step size. Song et al. \cite{szz} analyzed the pth moment asymptotical ultimate boundedness of PSDEs with time-varying coefficients.  Hu et al. \cite{HRH} established the existence and uniqueness for a class of  PSDEs. For more details, we refer the reader to \cite{ Z0,  LD, RHX1,M1,YYX} and references therein.

 Recently, authors \cite{why} have developed the fundamental theories for pantograph stochastic functional differential equations (PSFDEs) described by the following stochastic differential equation:
 \begin{align*}
\begin{cases}
 &\dif x(t)= f( x_t, t)\dif t +g(  x_t, t)\dif  B(t), t\in [0, \infty),\\
 & x(0)=x_{0},
 \end{cases}
 \end{align*}
 where   $x_{t}: =\{x(\theta t), \, \underline{\theta}\leq \theta\leq 1\}$ is a segment process and $0<\underline{\theta}<1$ is a fixed constant  and (more details can be seen in Section 2 below). One can see that  {\it the PSFDEs differs markedly from PSDEs, since the current state of PSFDEs depends on a past segment of its solution while the current state of PSDEs depends only on a past point of its solution.}
The existence and uniqueness, exponential stability and polynomial stability of solutions to \eqref{eq41} are investigated in \cite{why}. In the present paper, we shall investigate the numerical solutions of PSFDEs.
Our main contribution are as follows:
\begin{itemize}

\item[$\bullet$] We are the first to define  the numerical solution for PSFDEs. An implementable scheme is proposed such that  the numerical solution converges strongly to the analytical solution in finite time interval.

\item[$\bullet$] The numerical solutions preserve the exponential stability  and almost surely polynomial stability of the analytical solution under the certain conditions.

\end{itemize}

We close this part by giving our organization  in this article. In Section 2, we introduce some necessary notations. In Section 3, we show the strong convergence of the numerical solutions. The exponential stability  and almost surely polynomial stability  of numerical solutions are discussed in Section 4.  Several  examples are presented to illustrate the theories.

 \section{Notations and Preliminaries}
Throughout this paper,  Let  $(\Omega, \sF, \{\sF_{t}\}_{t\geq 0}, P)$ be  a complete probability space satisfying the usual conditions(i.e., it is increasing and right continuous with $\sF_{0}$ contains all $P$-null sets) and taking along
 a standard  $d$-Brownian motion process  $B(t).$  For $x, y \in \mR^{n},$ we use $|x |$  to denote the Euclidean norm of
$x,$ and use  $\langle x, y\rangle$ or $x^{T}y$ to denote the Euclidean inner product. If  $A$ is a matrix, $A^{T}$  is the transpose  of $A$  and   $|A |$ represents  $\sqrt{\mathrm{Tr} (AA^{T}).}$   Let $\lfloor a \rfloor$ be the integer parts of $a.$  Moreover, for $0<\underline{\theta} <1,$ denote by $\mathscr{C} :=\mathscr{C}([\underline{\theta}, 1]; \mR^{n})$ the family of all continuous $\mR^{n}-$valued functions $\varphi$  defined on $[\underline{\theta}, 1]$ with the norm $\| \varphi\|=\sup_{\underline{\theta} \leq t\leq 1}|\varphi(t)|.$

 Consider the following equations:
\begin{equation}\label{eq41}
\begin{cases}
 &\dif x(t)= f( x_t, t)\dif t +g(  x_t, t)\dif  B(t), t\in [0, \infty),\\
 & x(0)=x_{0},
 \end{cases}
 \end{equation}
where $x_{t}=\{x(\theta t), \underline{\theta}\leq \theta \leq 1)\}$ and $0< \underline{\theta}<1$ is a fixed constant. For convenience, we assume $f(0, t )=g(0,t)=0.$
 In \cite{why}, we have showed that the analytic solution to (\ref{eq41}) is exponential stable and  almost surely polynomial stable under some conditions.  In this paper, we will prove that the Euler-Maruyama( EM)  method can inherit  exponential stability  and almost surely polynomial stability of the analytical solution under the certain conditions, with the help of semi-martingale convergence theorem. Before giving our main results, we cite  the semi-martingale convergence theorem as a lemma which can be found in \cite{LR10}.

\bl
Let $\xi_{1}(i), \xi_{2}(i)$ be two sequences of nonnegative random variables such that $\xi_{1}(i), \xi_{2}(i)$ are $\sF_{i-1}-$measurable for $i=1,2,\cdots, $  with $\xi_{1}(0)=0, \xi_{2}(0)=0,$ a.s. Let  M(i) be a real-valued local martingale with $M(0)=0,$ a.s.,  and let $\xi$  be a nonnegative $\sF_{0}-$measurable random variable such that $\mE[\xi]< \infty.$ Set $\eta(i)=\xi+\xi_{1}(i)-\xi_{2}(i)+M(i), t\geq 0.$  If $\eta(i)$ is nonnegative , then we have the following results:
$$\{\lim_{i\rightarrow \infty}\xi_{1}(i)<\infty\}\subset \{\lim_{t\rightarrow \infty}\xi_{2}(i)<\infty\}\cap \{\lim_{t\rightarrow \infty}\eta(i)<\infty\},\,a.s.,$$
where $C\subset D,\, a.s.$  means $P(C\cup D^{c})=0.$ In particular, if $\lim_{t\rightarrow \infty}\xi_{1}(i)<\infty,\,a.s.,$ then, with probability one,
$$\lim_{i\rightarrow \infty}\xi_{2}(i)<\infty, \,  \lim_{i\rightarrow \infty}\eta(i)<\infty, -\infty <\lim_{i\rightarrow \infty}M(i)<\infty  ,\,a.s.$$
\el

\section{The EM method and strong convergence}

 Choose a step size $\Delta\in (0,1)$ and define the discrete EM approximate solution $y(k)=y(k\Delta)\approx x(k\Delta)$ by setting $y(0)=x_{0}, y_{0}=x_{0}$  and forming
\begin{align}\label{neq}
y(k+1)=y(k)+f(y_k, k\Delta)\Delta+g(y_k,  k\Delta)\Delta B(k), \, k=0, 1, 2, \ldots
\end{align}
where $y(k)=y(k\Delta), \Delta B(k)=B((k+1)\Delta)-B(k\Delta),$  and $y_{k} $  is  a $\mathscr{C}([\underline{\theta},1]; \mR^{n})-$valued random variable defined as follows:

\begin{equation}\label{i}
y_{k}(u)=\begin{cases}
&y((\lfloor  k\underline{\theta}\rfloor +i) \Delta)+\frac{u-(\lfloor  k\underline{\theta}\rfloor +i)\Delta}{\Delta}[y((\lfloor  k\underline{\theta}\rfloor +i+1) \Delta)-y((\lfloor  k\underline{\theta}\rfloor +i) \Delta)],\\
 &\quad \quad\quad \mbox{for}\, (\lfloor  k\underline{\theta}\rfloor +i)\Delta \leq u \leq (\lfloor  k\underline{\theta}\rfloor +i+1)\Delta , i=1,\cdots, k-\lfloor  k\underline{\theta}\rfloor-1;\\
 &y(\lfloor  k\underline{\theta}\rfloor  \Delta)+\frac{u-\lfloor  k\underline{\theta}\rfloor \Delta}{\Delta}[y((\lfloor  k\underline{\theta}\rfloor +1) \Delta)-y(\lfloor  k\underline{\theta}\rfloor  \Delta)],\\
 &\quad \quad\quad \mbox{for}\,   k\underline{\theta} \Delta \leq u \leq (\lfloor  k\underline{\theta}\rfloor +1)\Delta.
 \end{cases}
\end{equation}

\eqref{i} can be  rewritten as
\begin{equation}\label{ii}
y_{k}(u)=\begin{cases}
&\frac{\Delta -u+(\lfloor  k\underline{\theta}\rfloor +i)\Delta}{\Delta}y((\lfloor  k\underline{\theta}\rfloor +i) \Delta)+\frac{u-(\lfloor  k\underline{\theta}\rfloor +i)\Delta}{\Delta}y((\lfloor  k\underline{\theta}\rfloor +i+1) \Delta),\\
 &\quad \quad\quad \mbox{for}\, (\lfloor  k\underline{\theta}\rfloor +i)\Delta \leq u \leq (\lfloor  k\underline{\theta}\rfloor +i+1)\Delta , i=1,\cdots, k-\lfloor  k\underline{\theta}\rfloor-1;\\
 &\frac{\Delta-u+\lfloor  k\underline{\theta}\rfloor \Delta}{\Delta}y(\lfloor  k\underline{\theta}\rfloor  \Delta)+\frac{u-\lfloor  k\underline{\theta}\rfloor \Delta}{\Delta}y((\lfloor  k\underline{\theta}\rfloor +1) \Delta),\\
 &\quad \quad\quad \mbox{for}\,   k\underline{\theta} \Delta \leq u \leq (\lfloor  k\underline{\theta}\rfloor +1)\Delta.
 \end{cases}
\end{equation}
Thus,
\begin{align}\label{iii}
|y_{k}(u)|\leq |y((\lfloor  k\underline{\theta}\rfloor +i) \Delta)|\vee |y((\lfloor  k\underline{\theta}\rfloor +i+1) \Delta)|, (\lfloor  k\underline{\theta}\rfloor +i)\Delta \leq u \leq (\lfloor  k\underline{\theta}\rfloor +i+1)\Delta.
\end{align}
From \eqref{iii}, we have
\begin{align}\label{iv}
\|y_{k}\|\leq \sup_{\theta \in [\underline{\theta}, 1]}|y(\lfloor k\theta\rfloor \Delta)|.
\end{align}
In order  to analyze the continuous-time approximates , for $t\in [0, T],$ we define
\begin{align*}
z(t)= \sum^{\infty}_{k=0}y(k )1_{[k\Delta,(k+1)\Delta) }(t),\,
\bar{t}=\sum^{\infty}_{k=0} k\Delta 1_{[k\Delta,(k+1)\Delta) }(t),
\end{align*}
$$\bar{z}(\theta,t)= \sum^{\infty}_{k=0}y(\lfloor \theta k \rfloor )1_{[k\Delta,(k+1)\Delta) }(t).$$
We also define $\bar z_t$  a segment process on $\mathscr{C}([\underline{\theta},1]; \mR^{n})$  as the following
\begin{align*}
\bar{z}_{t}(\theta)=\sum^{\infty}_{k=0}y_{k}(\theta k\Delta)1_{[k\Delta,(k+1)\Delta) }(t), \, \theta \in [\underline{\theta},1].
\end{align*}
The continuous-time approximation $\{Y(t), t\geq 0\}$ is defined as $Y(0)=x_{0}$ and
\begin{align}\label{c}
Y(t)=x_{0}+\int^{t}_{0}f(\bar{z}_{s}, \bar{s})\dif s+\int^{t}_{0}g(\bar{z}_{s},\bar{s})\dif B(s).
\end{align}
For the future use, we make the following assumption.
\begin{enumerate}
\item[(H)]  There exists  positive constant $K$  such that  for any $\varphi \in C([\underline{\theta}, 1]; \mR^{n}) $
    \begin{align*}
    &|f( \varphi, t)-f( \phi, t) |^{2}\vee|g( \varphi, t) -g( \phi, t)|^{2}
    \leq K\|\varphi(\theta)-\phi(\theta)\|^{2}.
    \end{align*}
\end{enumerate}
\bl\label{L1}
Assume (H).
Then  there exists a positive constant $K_{1}$ which is only dependent  on  $K, x_{0},T, p\geq 2$ but independent of $\Delta$ such that
$$\mE\sup_{0\leq t\leq T}|Y(t)|^{p}\vee \mE\sup_{0\leq t\leq T}|x(t)|^{p}\le K_{1}.$$
\el
This  lemma can be proved by  similar method to that of  Lemma 3.2 in  Mao \cite{MXX}, we omit it here.
Now, we  state the  second lemma in this section.
\bl\label{L2}
Assume (H). Then   there exists a constant $K_{2}$ which is only dependent  on  $ K,K_{1}, x_{0},T$ but independent of $\Delta$ such that
$$\mE\sup_{0\leq t\leq T}|Y(t)-z(t)|^{2}\le K_{2}\Delta, \,  \mE\sup_{0\leq t\leq T}\|Y_{t}-\bar{z}_{t}\|^{2} \le K_{2}\Delta. $$
\el
\begin{proof}
For any $t\in [0,T],$ let $k=\lfloor \frac{t}{\Delta}\rfloor, m=\lfloor \frac{T}{\Delta}\rfloor,$  obviously, $k\leq m.$ Then,     $t\in [k\Delta, (k+1)\Delta), $ Thus, we have
\begin{align*}
Y(t)-z(t)=f(\bar{z}_{k\Delta},  k\Delta)(t-k\Delta)+g(\bar{z}_{_{k\Delta}},  k\Delta)(B(t)-B(k\Delta)).
    \end{align*}
According to assumption (H),  we obtain
\begin{align*}
|Y(t)-z(t)|^{2}&\leq 2K\|\bar{z}_{k\Delta}\|^{2}\Delta^{2}+2|g(\bar{z}_{_{k\Delta}},  k\Delta)(B(t)-B(k\Delta))|^{2}.
 \end{align*}
Thus,
\begin{align}\label{3.4}
&\mE\sup_{k\Delta\leq t\leq (k+1)\Delta}|Y(t)-z(t)|^{2}\leq 2K\mE[\sup_{k\Delta\leq t\leq (k+1)\Delta}\|\bar{z}_{t}\|^{2}]\Delta^{2}\no\\
&\quad \quad\quad\quad \quad\quad\quad \quad\quad\quad\quad+2\mE[\sup_{k\Delta\leq t\leq (k+1)\Delta}\quad|g(\bar{z}_{_{k\Delta}},  k\Delta)(B(t)-B(k\Delta))|^{2}]\no\\
&\leq 2K\mE[\sup_{k\Delta\leq t\leq (k+1)\Delta}\|\bar{z}_{t}\|^{2}]\Delta^{2}\no\\
&\quad+2(\mE\sup_{k\Delta\leq t\leq (k+1)\Delta}|g(\bar{z}_{_{k\Delta}},  k\Delta)|^{4})^{\frac{1}{2}}(\mE\sup_{k\Delta\leq t\leq (k+1)\Delta}|B(t)-B(k\Delta)|^{4})^{\frac{1}{2}}\no\\
&\leq 2KK_{1,2}\Delta^{2} +2KK^{\frac{1}{2}}_{1,4}(\mE\sup_{k\Delta\leq t\leq (k+1)\Delta}|B(t)-B(k\Delta)|^{4})^{\frac{1}{2}}\no\\
& \leq 2KK_{1,2}\Delta^{2} +2d^{\frac{1}{2}}KK^{\frac{1}{2}}_{1,4}\bigg(\bigg.\mE\sup_{t\in [k\Delta, (k+1)\Delta \wedge T]}|B_{j}(t)-B_{j}(k\Delta)|^{4}\bigg)\bigg.^{\frac{1}{2}}
\end{align}
By Doob's martingale inequality,  we have
\begin{align}\label{3.5}
\mE\sup_{k\Delta\leq t\leq (k+1)\Delta}|B_{j}(t)-B_{j}(k\Delta)|^{4}
\leq \frac{256}{27}\Delta^{2}.
\end{align}
By \eqref{3.4} and \eqref{3.5}, we arrive at
\begin{align}
\mE&\sup_{k\Delta\leq t\leq (k+1)\Delta}|Y(t)-z(t)|^{2}
\leq 2KK_{1,2}\Delta^{2} +2d^{\frac{3}{2}}KK^{\frac{1}{2}}_{1,4}\frac{16}{3\sqrt{3}} \Delta\no\\
&\leq \bigg[\bigg.2KK_{1,2} +d^{\frac{3}{2}}K^{2}K_{1,4}\frac{32}{3\sqrt{3}} \bigg]\bigg.\Delta.
\end{align}
Next, we prove the second result.
\begin{align}\label{3.7}
\mE&\sup_{k\Delta\leq t\leq (k+1)\Delta}\sup_{\underline{\theta}\leq \theta \leq 1}|Y(\theta t)-\bar{z}_{t}(\theta)|^{2}\no\\
&\leq  3\bigg[\bigg.2KK_{1,2} +d^{\frac{3}{2}}K^{2}K_{1,4}\frac{32}{3\sqrt{3}} \bigg]\bigg.\Delta+ 3\mE\sup_{t\in [k\Delta, (k+1)\Delta \wedge T]}\sup_{\underline{\theta}\leq \theta \leq 1}|z(\theta t)-\bar{z}(t,\theta)|^{2}\no\\
& \quad +3\mE\sup_{t\in [k\Delta, (k+1)\Delta \wedge T]}\sup_{\underline{\theta}\leq \theta \leq 1}|\bar{z}( t,\theta )-\bar{z}_{t}(\theta)|^{2}=:I_{1}+I_{2}+I_{3}.
\end{align}
 It is obvious that $\theta t\in [k^{\theta}_{t}\Delta, (k^{\theta}_{t}+1)\Delta],$ where $k^{\theta}_{t}=\lfloor\frac{\theta t}{\Delta}\rfloor.$  Thus, for $k=\lfloor\frac{t}{\Delta}\rfloor,$ we get
 \begin{align}\label{++}
 |z(\theta t)-\bar{z}(\theta, t)|=|y(k^{\theta}_{t}\Delta)-y(\lfloor \theta k\rfloor\Delta)|.
 \end{align}
Since $\frac{\theta t}{\Delta}\in [\theta k, \theta(k+1)),$  it can be seen that $\lfloor \theta k\rfloor\Delta\leq k^{\theta}_{t}\leq \lfloor \theta (k+1)\rfloor\leq \lfloor \theta k\rfloor +1.$ Then, $k^{\theta}_{t}- \lfloor \theta k\rfloor\leq 1.$   By \eqref{++},  one has
\begin{align}\label{3.8}
I_{2}=\mE&\sup_{t\in [k\Delta, (k+1)\Delta \wedge T]}\sup_{\underline{\theta}\leq \theta \leq 1}|z(\theta t)-\bar{z}(t,\theta)|^{2}
\leq  \mE\sup_{0\leq l\Delta\leq T}|y((l+1)\Delta\wedge T)-y(l\Delta)|^{2}.
\end{align}
Next, we calculate $I_{3}.$  From \eqref{i}, we have
\begin{align}\label{+++}
I_{3}=\mE\sup_{t\in [k\Delta, (k+1)\Delta \wedge T]}\sup_{\underline{\theta}\leq \theta \leq 1}|\bar{z}( t,\theta )-\bar{z}_{t}(\theta)|^{2}
\leq  \mE\sup_{0\leq l\Delta\leq T}|y((l+1)\Delta\wedge T)-y(l\Delta)|^{2}.
\end{align}
By using similar procedure as used method in the proof of first result, one can see that
\begin{align}\label{3.9}
\mE\sup_{0\leq l\Delta\leq T}|y((l+1)\Delta\wedge T)-y(l\Delta)|^{2}\leq K_{3}\Delta, \end{align}
where $K_{3}$ is only dependent  on  $K,K_{1,p}, x_{0},T$ but independent of $\Delta$.
\eqref{3.7}, \eqref{++},  \eqref{+++}, \eqref{3.8} and \eqref{3.9} lead to the second result.
\end{proof}

The following result reveals that the numerical solutions converge to the true solution.
\bt
Assume (H).   It holds that
$$\mE\sup_{0\leq t\leq T}|x(t)-Y(t)|^{2}\leq K_{3}\Delta,$$
where $K_{4}$ is only dependent  on  $K, K_{1,p},K_{2},x_{0},T$ but independent of $\Delta.$
\et
\begin{proof}
By (H), Lemmas \ref{L1} and \ref{L2}, we compute
\begin{align*}
&\mE\sup_{0\leq s\leq t}|x(s)-Y(s)|^{2}\\
&\leq 2T\mE\int^{t}_{0}|f(x_{s}, s)-f(\bar{z}_{s}, s)|^{2}\dif s
 +2\mE\sup_{0\leq s\leq t}\bigg|\bigg.\int^{s}_{0}g(x_{s},s)-g(\bar{z}_{s}, s)\dif B(s)   \bigg|\bigg.^{2}\\
&\leq 2T\mE\int^{t}_{0}|f(x_{s}, s)-f(\bar{z}_{s}, s)|^{2}\dif s +8\mE\int^{t}_{0}|g(x_{s}, s)-g(\bar{z}_{s}, s)|^{2}\dif s\\
&\leq 4TK\mE\int^{t}_{0}\|x_{s}-Y_{s}\|^{2}\dif s  +4TK\mE\int^{t}_{0}\|Y_{s}-\bar{z}_{s}\|^{2}\dif s\\
&+16K\mE\int^{t}_{0}\|x_{s}-Y_{s}\|^{2}\dif s  +16K\mE\int^{t}_{0}\|Y_{s}-\bar{z}_{s}\|^{2}\dif s\\
&\leq 4TK\mE\int^{t}_{0}\mE\sup_{0\leq u\leq s}|x(u)-Y(u)|^{2}\dif s  +4TK\mE\int^{t}_{0}\|Y_{s}-\bar{z}_{s}\|^{2}\dif s\\
&+16K\int^{t}_{0}\mE\sup_{0\leq u\leq s}|x(u)-Y(u)|^{2}\dif s  +16K\mE\int^{t}_{0}\|Y_{s}-\bar{z}_{s}\|^{2}\dif s\\
&\leq C\Delta
+(4TK+16K)\int^{t}_{0}\mE\sup_{0\leq u\leq s}|x(u)-Y(u)|^{2}\dif s,
\end{align*}
where $C$ is a positive constant independent of stepsize $\Delta$. Gronwall's inequality leads to required result.
\end{proof}

\section{ Stability of numerical solutions}

In the section, we shall investigate the exponential stability and polynomial stability for the numerical solutions.

\subsection{Exponential stability of numerical solution}

We need the  following assumptions.
\begin{enumerate}
\item[(H1)]  For any $\varphi \in C([\underline{\theta}, 1]; \mR^{n}),$  there exists a probability measures $\nu$ on $[\underline{\theta}, 1]$  with positive constants $\lambda_{1}, \lambda_{2} $ such that
  \begin{align}
&2\langle \varphi(1)-\phi(1), f(\varphi,t)-f(\phi,t)\rangle+|g(\varphi,t)-g(\phi,t)|^{2}\no\\
&\leq -\lambda_{1}|\varphi(1)-\phi(1)|^{2}+\lambda_{2}\int^{1}_{\underline{\theta}}e^{-\beta t}|\varphi(\theta)-\phi(\theta)|^{2}\dif \nu(\theta).
 \end{align}
\end{enumerate}
\begin{enumerate}
\item[(H2)] For any $\varphi, \phi \in C([\underline{\theta}, 1]; \mR^{n}),$   there exists a probability measures $\nu$ on $[\underline{\theta}, 1]$  with positive constants $\lambda_{3}$ and  $\lambda_{4}$ such that
    \begin{align*}
    &|f( \varphi, t)-f( \phi, t) |^{2}\vee|g( \varphi, t)-g(\phi, t) |^{2}\\
    &\leq \lambda_{3}|\varphi(1)-\phi(1)|^{2}+\lambda_{4}\int^{1}_{\underline{\theta}} e^{-\beta t}|\varphi(\theta)-\phi(\theta)|^{2}\dif \nu(\theta),
    \end{align*}
where  $\beta$ is a constant satisfying $0<\frac{1-\beta}{\underline{\theta}}<1.$
\end{enumerate}

We can see that (H2) implies (H).

\bt \label{YT1}
Assume  $\mathrm{(H1)}$ and  $\mathrm{(H2)}$.
If   the following conditions hold:
\begin{enumerate}
\item[i)] there exist some positive constants $\bar{C}, \alpha_{0}$ satisfying   $1<\bar{C}\leq e^{\alpha_{0}}$  and a  sufficiently small constant  $ \lambda_{0}> 0$ such that
\begin{align*}
H(\bar{C}, \lambda_{0})= -\lambda_{1}+\alpha_{0}+ 2\lambda_{2}\bigg(\bigg.\bigg\lfloor\bigg.\frac{1}{\underline{\theta}}\bigg\rfloor\bigg. +1\bigg)\bigg.\bar{C}^{\frac{\lambda_{0}}{\underline{\theta}}}
 +\lambda_{3}\lambda_{0} +2\lambda_{4}\bigg(\bigg.\bigg\lfloor\bigg.\frac{1}{\underline{\theta}}\bigg\rfloor\bigg. +1\bigg)\bigg.\bar{C}^{\frac{\lambda_{0}}{\underline{\theta}}}\lambda_{0}\leq 0.
\end{align*}
\item[ii)] $\Delta \in (0,\lambda_{0})$ is  small enough   satisfying the following inequality:
  $$\bigg[\bigg. \lambda_{2}\bigg(\bigg.\bigg\lfloor\bigg.\frac{1}{\underline{\theta}}\bigg\rfloor\bigg. +1\bigg)\bigg.\bar{C}^{\frac{\Delta}{\underline{\theta}}}+\lambda_{4}\bigg(\bigg.\bigg\lfloor\bigg.\frac{1}{\underline{\theta}}\bigg\rfloor\bigg. +1\bigg)\bigg.\bar{C}^{\frac{\Delta}{\underline{\theta}}}\Delta\bigg]\bigg.\bar{C}^{\Delta}\Delta\leq \frac{1}{2}.$$
\end{enumerate}   Then
the approximate solution $y(k)$ satisfies
  \begin{equation*}
  \begin{split}
&\limsup_{k\rightarrow \infty }\frac{1}{k}\log|y(k)|^{2}\leq -\alpha,\\
&\limsup_{k\rightarrow \infty }\frac{1}{k}\log\mE[|y(k)|^{2}]\leq -\alpha,
\end{split}
\end{equation*}
where $\alpha $ is a constant satisfying $e^{\alpha}=\bar{C}.$
\et

\begin{proof}
By virtue of  (\ref{neq}),  we have
\begin{equation}
\begin{split}
|y(k+1)|^{2}&\leq |y(k)|^{2}+ |f(y_k,k\Delta)|^{2}\Delta^{2}+2y^{T}(k)f(y_k,k\Delta)\Delta\\
&+|g(y_k,k\Delta)|^{2}\Delta+|g(y_k,k\Delta)|^{2}((\Delta B(k))^{2}-\Delta)
\\
&+2y^{T}(k)g(y_k,k\Delta)\Delta B(k)+2f(y_k,k\Delta)g(y_k,k\Delta)\Delta B(k)\Delta\\
&\leq |y(k)|^{2}+ |f(y_k,k\Delta)|^{2}\Delta^{2}+2y^{T}(k)f(y_k,k\Delta)\Delta\\
&+|g(y_k,k\Delta)|^{2}\Delta+M(k),
\end{split}
\end{equation}
where
\begin{align*}
M(k)&=|g(y_k,k\Delta)|^{2}((\Delta B(k))^{2}-\Delta)+2y^{T}(k)g(y_k,k\Delta)\Delta B(k)\\
&+2 f(y_k,k\Delta)g(y_k,k\Delta)\Delta B(k)\Delta.
\end{align*}
By (H1) and  (H2),  one can see that
\begin{align}\label{3.12}
&|y(k+1)|^{2}-|y(k)|^{2}\leq\bigg(\bigg. -\lambda_{1}|y(k)|^{2}+ \lambda_{2}\int^{1}_{\underline{\theta}}e^{-\beta k\Delta}|y_{k}( k\theta  \Delta)|^{2}\dif \nu(\theta)\bigg)\bigg. \Delta \no\\
 &+\bigg(\bigg. \lambda_{3}|y(k)|^{2}+ \lambda_{4}\int^{1}_{\underline{\theta}}e^{-\beta k\Delta}|y_{k}( k\theta\Delta)|^{2}\dif \nu(\theta)\bigg)\bigg. \Delta^{2} +M(k).
\end{align}
Multiplying $C^{(j+1)\Delta}$ on both sides of the inequality \eqref{3.12} yields that
\begin{align}\label{3.13}
&C^{(j+1)\Delta}|y(j+1)|^{2}-C^{j\Delta}|y(j)|^{2}\no\\
&=C^{(j+1)\Delta}(1-\frac{1}{C^{\Delta}})|y(j)|^{2}\no\\
&+\bigg(\bigg. -\lambda_{1}C^{(j+1)\Delta}|y(j)|^{2}+ \lambda_{2}\int^{1}_{\underline{\theta}}C^{(j+1)\Delta}e^{-\beta j\Delta}|y_{j}( j\theta  \Delta)|^{2}\dif \nu(\theta)\bigg)\bigg. \Delta \no\\
&+\bigg(\bigg. \lambda_{3}C^{(j+1)\Delta}|y(j)|^{2}+ \lambda_{4}\int^{1}_{\underline{\theta}}C^{(j+1)\Delta}e^{-\beta j\Delta}|y_{j}( j\theta  \Delta)|^{2}\dif \nu(\theta)\bigg)\bigg. \Delta^{2} +M(j),
\end{align}
where $C$ is a constant satisfying $1<C\leq e^{\alpha_{0}}.$  Since $1-C^{-\Delta}<\alpha_{0}\Delta,$
summing both sides of \eqref{3.13} from $j=0$ to $j=k-1,$  we obtain
\begin{align}\label{3.14}
&C^{k\Delta}|y(k)|^{2}\no\\
&=x_{0}+\sum^{k-1}_{j=0}C^{(j+1)\Delta}(1-\frac{1}{C^{\Delta}})|y(j)|^{2}+\bigg(\bigg. -\lambda_{1}\sum^{k-1}_{j=0}C^{(j+1)\Delta}|y(j)|^{2}\no\\
&+ \lambda_{2}\sum^{k-1}_{j=0}\int^{1}_{\underline{\theta}}C^{(j+1)\Delta}e^{-\beta j\Delta}|y_{j}( j\theta  \Delta)|^{2}\dif \nu(\theta)\bigg)\bigg. \Delta \no\\
&+\bigg(\bigg. \lambda_{3}\sum^{k-1}_{j=0}C^{(j+1)\Delta}|y(j)|^{2}+ \lambda_{4}\int^{1}_{\underline{\theta}}\sum^{k-1}_{j=0}C^{(j+1)\Delta}e^{-\beta j\Delta}|y_{j}( j\theta  \Delta)|^{2}\dif \nu(\theta)\bigg)\bigg. \Delta^{2} +\sum^{k-1}_{j=0}M(j),
\end{align}
where $\sum^{k-1}_{j=0}M(j)$ is a martingale.  Firstly,  we compute
\begin{align}
\sum^{k-1}_{j=0}\int^{1}_{\underline{\theta}}C^{(j+1)\Delta}e^{-\beta j\Delta}|y_{j}( j\theta  \Delta)|^{2}\dif\nu(\theta) =\int^{1}_{\underline{\theta}}\sum^{k-1}_{j=0}C^{(j+1)\Delta}e^{-\beta j\Delta}|y_{j}( j\theta  \Delta)|^{2}\dif \nu(\theta).
\end{align}
It is not difficult we can see that
\begin{align}\label{3.16}
\sum^{k-1}_{j=0}&C^{(j+1)\Delta}e^{-\beta j\Delta}|y(\lfloor j\theta \rfloor )|^{2}\no\\
&=|y(0) |\sum_{0\leq j \leq k-1:\lfloor \theta j\rfloor =0}C^{(j+1)\Delta}e^{-\beta j\Delta}+|y(1) |\sum_{0\leq j \leq k-1:\lfloor \theta j\rfloor =1}C^{(j+1)\Delta}e^{-\beta j\Delta}\no\\
&+\cdots +|y(\lfloor \theta (k-1)\rfloor) |\sum_{0\leq j \leq k-1:\lfloor \theta j\rfloor =\lfloor \theta (k-1)\rfloor}C^{(j+1)\Delta}e^{-\beta j\Delta},
\end{align}
and
\begin{align}\label{3.16+}
\sum^{k-1}_{j=0}&C^{(j+1)\Delta}e^{-\beta j\Delta}|y(\lfloor j\theta \rfloor+1)|^{2}\no\\
&=|y(1) |\sum_{0\leq j \leq k-1:\lfloor \theta j\rfloor =0}C^{(j+1)\Delta}e^{-\beta j\Delta}+|y(2) |\sum_{0\leq j \leq k-1:\lfloor \theta j\rfloor =1}C^{(j+1)\Delta}e^{-\beta j\Delta}\no\\
&+\cdots +|y(\lfloor \theta (k-1)\rfloor+1) |\sum_{0\leq j \leq k-1:\lfloor \theta j\rfloor =\lfloor \theta (k-1)\rfloor}C^{(j+1)\Delta}e^{-\beta j\Delta}.
\end{align}

Additionally, noting  that $\lfloor \theta j\rfloor  =i\Leftrightarrow \frac{i}{\theta}\leq j <\frac{i+1}{\theta},$ for any $i=0,1,\cdots, \lfloor \theta (k-1)\rfloor,$ which implies
\begin{align*}
\begin{cases}
\frac{i}{\theta}\leq j \leq\frac{i+1}{\theta}-1, \frac{1}{\theta}\in \mN,\\
\lfloor\frac{i}{\theta}\rfloor +1\leq j \leq \lfloor\frac{i+1}{\theta}\rfloor\leq \lfloor\frac{i}{\theta}\rfloor + \lfloor\frac{1}{\theta}\rfloor +1, \frac{1}{\theta} \notin \mN.
\end{cases}
\end{align*}
Then,  it yields that the number of those $j\in \{0,1,2,\cdots, k-1\}$  such that $\lfloor \theta j\rfloor =i,$ for some $i\in \{0,1,2,\cdots, \lfloor \theta (k-1)\rfloor\}$, is at most $\lfloor\frac{1}{\theta}\rfloor +1.$  Moreover, the greatest $j$  for which $\lfloor \theta j\rfloor=i$ is less than $\frac{i+1}{\theta}$  and greater that $\frac{i}{\theta}.$  By \eqref{3.16}, we derive that
\begin{align}\label{1212}
&\sum_{j=0}^{k-1}C^{(j+1)\Delta}e^{-\beta j\Delta}|y(\lfloor \theta j\rfloor)|^{2}\no\\
&\leq \bigg(\bigg.\bigg\lfloor\bigg.\frac{1}{\theta}\bigg\rfloor\bigg. +1\bigg)\bigg.\sum^{\lfloor\theta (k-1)\rfloor}_{j=0}C^{(\frac{j+1}{\theta}+1)\Delta}e^{-\beta \frac{j}{\theta}\Delta}|y(j)|^{2} \no\\
&\leq  \bigg(\bigg.\bigg\lfloor\bigg.\frac{1}{\theta}\bigg\rfloor\bigg. +1\bigg)\bigg.\sum^{\lfloor\theta (k-1)\rfloor}_{j=0}C^{(\frac{1-\beta}{\theta}j\Delta+\frac{\Delta}{\theta}+\Delta)}|y(j)|^{2} \no\\
& \leq \bigg(\bigg.\bigg\lfloor\bigg.\frac{1}{\theta}\bigg\rfloor\bigg. +1\bigg)\bigg.C^{\frac{\Delta}{\theta}}\sum^{\lfloor\theta (k-1)\rfloor}_{j=0}C^{(j+1)\Delta}|y(j)|^{2}\no\\
&\leq \bigg(\bigg.\bigg\lfloor\bigg.\frac{1}{\underline{\theta}}\bigg\rfloor\bigg. +1\bigg)\bigg.C^{\frac{\Delta}{\underline{\theta}}}\sum^{k-1}_{j=0}C^{(j+1)\Delta}|y(j)|^{2},
\end{align}
and
\begin{align}\label{1212+}
&\sum_{j=0}^{k-1}C^{(j+1)\Delta}e^{-\beta j\Delta}|y(\lfloor \theta j\rfloor+1)|^{2}\no\\
&\leq \bigg(\bigg.\bigg\lfloor\bigg.\frac{1}{\theta}\bigg\rfloor\bigg. +1\bigg)\bigg.\sum^{\lfloor\theta (k-1)\rfloor}_{j=0}C^{(\frac{j+1}{\theta}+1)\Delta}e^{-\beta \frac{j}{\theta}\Delta}|y(j+1)|^{2} \no\\
&\leq  \bigg(\bigg.\bigg\lfloor\bigg.\frac{1}{\theta}\bigg\rfloor\bigg. +1\bigg)\bigg.\sum^{\lfloor\theta (k-1)\rfloor}_{j=0}C^{(\frac{1-\beta}{\theta}j\Delta+\frac{\Delta}{\theta}+\Delta)}|y(j+1)|^{2} \no\\
& \leq \bigg(\bigg.\bigg\lfloor\bigg.\frac{1}{\theta}\bigg\rfloor\bigg. +1\bigg)\bigg.C^{\frac{\Delta}{\theta}}\sum^{\lfloor\theta (k-1)\rfloor}_{j=0}C^{(j+1)\Delta}|y(j+1)|^{2}\no\\
&\leq \bigg(\bigg.\bigg\lfloor\bigg.\frac{1}{\underline{\theta}}\bigg\rfloor\bigg. +1\bigg)\bigg.C^{\frac{\Delta}{\underline{\theta}}}\sum^{k}_{j=1}C^{j\Delta}|y(j)|^{2},
\end{align}
Combining with \eqref{3.16}  and \eqref{3.16+},  one has
 \begin{align}
&C^{k\Delta}|y(k)|^{2}\no\\
&\leq x_{0}+\sum^{k-1}_{j=0}C^{(j+1)\Delta}\bigg(\bigg. 1-\frac{1}{C^{\Delta}}\bigg)\bigg.|y(j)|^{2}+\bigg[\bigg. -\lambda_{1}\sum^{k-1}_{j=0}C^{(j+1)\Delta}|y(j)|^{2}\no\\
&+ \lambda_{2}\bigg(\bigg.\bigg\lfloor\bigg.\frac{1}{\underline{\theta}}\bigg\rfloor\bigg. +1\bigg)\bigg.C^{\frac{\Delta}{\underline{\theta}}}\sum^{k-1}_{j=0}C^{(j+1)}|y(j)|^{2}\bigg]\bigg. \Delta \no\\
&+\bigg[\bigg. \lambda_{3}\sum^{k-1}_{j=0}C^{(j+1)\Delta}|y(j)|^{2}+ \lambda_{4}\bigg(\bigg.\bigg\lfloor\bigg.\frac{1}{\underline{\theta}}\bigg\rfloor\bigg. +1\bigg)\bigg.C^{\frac{\Delta}{\underline{\theta}}}\sum^{k-1}_{j=0}C^{(j+1)}|y(j)|^{2}\bigg]\bigg. \Delta^{2} \no\\
&+ \lambda_{2}\bigg(\bigg.\bigg\lfloor\bigg.\frac{1}{\underline{\theta}}\bigg\rfloor\bigg. +1\bigg)\bigg.C^{\frac{\Delta}{\underline{\theta}}}\sum^{k}_{j=1}C^{j\Delta}|y(j)|^{2} \Delta+ \lambda_{4}\bigg(\bigg.\bigg\lfloor\bigg.\frac{1}{\underline{\theta}}\bigg\rfloor\bigg. +1\bigg)\bigg.C^{\frac{\Delta}{\underline{\theta}}}\sum^{k}_{j=1}C^{j\Delta}|y(j)|^{2} \Delta^{2}\no\\
&+\sum^{k-1}_{j=0}M(j)\\
&\leq x_{0}+\bigg[\bigg. -\lambda_{1}+\alpha_{0}+ \lambda_{2}\bigg(\bigg.\bigg\lfloor\bigg.\frac{1}{\underline{\theta}}\bigg\rfloor\bigg. +1\bigg)\bigg.C^{\frac{\Delta}{\underline{\theta}}}\no\\
&\quad \quad \quad \quad\quad \quad\quad \quad\quad \quad\quad \quad+\lambda_{3}\Delta +\lambda_{4}\bigg(\bigg.\bigg\lfloor\bigg.\frac{1}{\underline{\theta}}\bigg\rfloor\bigg. +1\bigg)\bigg.C^{\frac{\Delta}{\underline{\theta}}}\Delta\bigg]\bigg.\sum^{k-1}_{j=0}C^{(j+1)\Delta}|y(j)|^{2} \Delta\no\\
&+\bigg[\bigg.  \lambda_{2}\bigg(\bigg.\bigg\lfloor\bigg.\frac{1}{\underline{\theta}}\bigg\rfloor\bigg. +1\bigg)\bigg.C^{\frac{\Delta}{\underline{\theta}}}+\lambda_{4}\bigg(\bigg.\bigg\lfloor\bigg.\frac{1}{\underline{\theta}}\bigg\rfloor\bigg. +1\bigg)\bigg.C^{\frac{\Delta}{\underline{\theta}}}\Delta\bigg]\bigg.\sum^{k}_{j=1}C^{j\Delta}|y(j)|^{2} \Delta+\sum^{k-1}_{j=0}M(j)\no\\
&\leq x_{0}+\bigg[\bigg. -\lambda_{1}+\alpha_{0}+ 2\lambda_{2}\bigg(\bigg.\bigg\lfloor\bigg.\frac{1}{\underline{\theta}}\bigg\rfloor\bigg. +1\bigg)\bigg.C^{\frac{\Delta}{\underline{\theta}}}\no\\
&\quad \quad \quad \quad\quad +\lambda_{3}\Delta +2\lambda_{4}\bigg(\bigg.\bigg\lfloor\bigg.\frac{1}{\underline{\theta}}\bigg\rfloor\bigg. +1\bigg)\bigg.C^{\frac{\Delta}{\underline{\theta}}}\Delta\bigg]\bigg.\sum^{k-1}_{j=0}C^{(j+1)\Delta}|y(j)|^{2} \Delta\no\\
&+\bigg[\bigg.  \lambda_{2}\bigg(\bigg.\bigg\lfloor\bigg.\frac{1}{\underline{\theta}}\bigg\rfloor\bigg. +1\bigg)\bigg.C^{\frac{\Delta}{\underline{\theta}}}+\lambda_{4}\bigg(\bigg.\bigg\lfloor\bigg.\frac{1}{\underline{\theta}}\bigg\rfloor\bigg. +1\bigg)\bigg.C^{\frac{\Delta}{\underline{\theta}}}\Delta\bigg]\bigg.C^{k\Delta}|y(k)|^{2} \Delta+\sum^{k-1}_{j=0}M(j).
\end{align}
Set
 \begin{align}
H(C, \Delta)= -\lambda_{1}+\alpha_{0}+ 2\lambda_{2}\bigg(\bigg.\bigg\lfloor\bigg.\frac{1}{\underline{\theta}}\bigg\rfloor\bigg. +1\bigg)\bigg.C^{\frac{\Delta}{\underline{\theta}}}
+\lambda_{3}\Delta +2\lambda_{4}\bigg(\bigg.\bigg\lfloor\bigg.\frac{1}{\underline{\theta}}\bigg\rfloor\bigg. +1\bigg)\bigg.C^{\frac{\Delta}{\underline{\theta}}}\Delta.
\end{align}
Then,
\begin{align}
H(\bar{C}, \lambda_{0})= -\lambda_{1}+\alpha_{0}+ 2\lambda_{2}\bigg(\bigg.\bigg\lfloor\bigg.\frac{1}{\underline{\theta}}\bigg\rfloor\bigg. +1\bigg)\bigg.\bar{C}^{\frac{\lambda_{0}}{\underline{\theta}}}
 +\lambda_{3}\lambda_{0} +2\lambda_{4}\bigg(\bigg.\bigg\lfloor\bigg.\frac{1}{\underline{\theta}}\bigg\rfloor\bigg. +1\bigg)\bigg.\bar{C}^{\frac{\lambda_{0}}{\underline{\theta}}}\lambda_{0}.
\end{align}
By condition i),  we have
$$H(\bar{C},\lambda_{0})\leq 0.$$
Using  condition ii)  and choosing a constant $\alpha>0$ with $e^{\eta}=\bar{C},$  one has
\begin{align*}
e^{\eta k\Delta}|y(k)|^{2}
\leq 2x_{0}+ 2\sum^{k-1}_{j=0}M(j).
\end{align*}
Since $\sum^{k-1}_{j=0}M(j)$ is a martingale, we get
$$\limsup_{k\rightarrow \infty }e^{\eta k\Delta}\mE|y(k)|^{2}
< \infty.$$
Furthermore, by Lemma 2.1 leads to
\begin{align}
\limsup_{k\rightarrow \infty }e^{\eta k\Delta}|y(k)|^{2}
< \infty.
\end{align}
The proof is therefore complete.
\end{proof}
\br
  Under the conditions of Theorem \ref{YT1},  Theorem 3.3 in \cite{why} has shown  that the analytical solution of \eqref{eq41} has the property of exponential stability. This means that the EM numerical solutions \eqref{neq} preserves the property of exponential stability of the analytical solution of \eqref{eq41}.
\er
We now give an example to explain Theorem \ref{YT1}.
\begin{exa} {\rm
Consider the following equation:
\begin{align}\label{EY1}
 &\dif x(t)= f( x_t, t)\dif t +g(  x_t, t)\dif  B(t), t\in [0, \infty)
 & x(0)=x_{0},
 \end{align}
where
\begin{align*}
f( \varphi, t)=
-1.1 \varphi(1)+0.04\int^{1}_{\frac{3}{4}}e^{-0.7t}|\varphi(\theta )|\dif \nu (\theta).
 \end{align*}
and
\begin{align*}
g( \varphi, t)=
0.2\int^{1}_{\frac{3}{4}}e^{-0.7t}|\varphi(\theta )|\dif \nu (\theta).
\end{align*}
Then,
\begin{align*}
&2\langle \varphi(1)-\phi(1), f(\varphi, t)-f(\phi,t)\rangle +|g(\varphi,t)-g(\phi,t)|^{2}\\
&=2\langle \varphi(1)-\phi(1), -1.1(\varphi(1)-\phi(1))+0.04\int^{1}_{\frac{3}{4}}e^{-0.7t}(\varphi(\theta )-\phi(\theta ))\dif \nu (\theta)\rangle \\
&+\bigg|\bigg.0.2\int^{1}_{\frac{3}{4}}e^{-0.7t}(\varphi(\theta t)-\phi(\theta t))\dif \nu (\theta)\bigg|\bigg.^{2}\\
& \leq -2.2|\varphi(1)-\phi(1)|^{2}+0.08(\varphi(1)-\phi(1))\int^{1}_{\frac{3}{4}}e^{-0.7t}(\varphi(\theta )-\phi(\theta ))\dif \nu (\theta)\\
&+\bigg|\bigg.0.2\int^{1}_{\frac{3}{4}}e^{-0.7t}(\varphi(\theta t)-\phi(\theta t))\dif \nu (\theta)\bigg|\bigg.^{2}\\
& \leq -2.16|\varphi(1)-\phi(1)|^{2}+0.08\int^{1}_{\frac{3}{4}}|e^{-0.7t}(\varphi(\theta )-\phi(\theta ))|^{2}\dif \nu (\theta),
\end{align*}
and
\begin{align*}
& |f(\varphi, t)-f(\phi,t)|^{2}\vee |g(\varphi, t)-g(\phi,t)|^{2}\\
& \leq 1.23|\varphi(1)-\phi(1)|^{2}+0.17\int^{1}_{\frac{3}{4}}e^{-0.7t}|(\varphi(\theta )-\phi(\theta ))|^{2}\dif \nu (\theta).
\end{align*}
We can find that
$$\lambda_{1}= 2.16,\lambda_{2}= 0.08, \lambda_{3}=1.23,\lambda_{4}= 0.17, \underline{\theta}=\frac{3}{4}.$$
Choosing $\bar{C}=1.1, \lambda_{0}=\frac{1}{300},\alpha_{0}=\frac{1}{10},$  it is  easily seen that the conditions i) and ii) of Thereem \ref{YT1} are satisfied.
Then, we conclude  that the numerical solutions of \eqref{EY1} are almost surely exponential stable,  and exponential stable in mean square.
}
\end{exa}

\subsection{ Almost sure polynomial stability of numerical solutions}
Next, we will study  polynomial stability of numerical solution to \eqref{neq}.    In this subsection, we assume that there exists a positive constant  $\overline{\theta}$ satisfying $\frac{1}{2}\vee \underline{\theta}<\bar{\theta}<1$.  We need the following assumptions.
\begin{enumerate}
\item[(H3)]  For any $\varphi \in C([\underline{\theta}, 1]; \mR^{n}),$  there exists a probability measure $\nu$ on $[\underline{\theta}, \overline{\theta}]$ and  positive constants $\lambda_{1}, \lambda_{2}$  such that
  \begin{align}
&2\langle \varphi(1)-\phi(1), f(\varphi,t)-f(\phi,t)\rangle+|g(\varphi,t)-g(\phi,t)|^{2}\no\\
&\leq -\lambda_{1}|\varphi(1)-\phi(1)|^{2}+\lambda_{2}\int^{\overline{\theta}}_{\underline{\theta}}|\varphi(\theta)-\phi(\theta)|^{2}\dif \nu(\theta).
 \end{align}
\end{enumerate}
\begin{enumerate}
\item[(H4)] For any $\varphi, \phi \in C([\underline{\theta}, 1]; \mR^{n}),$   there exists a probability measure $\nu$ on $[\underline{\theta}, \overline{\theta}]$ and two positive constants $\lambda_{3}, \lambda_{4}$ such that
    \begin{align*}
   & |f( \varphi, t)- f( \phi, t)|^{2}\vee|g( \varphi, t) -f( \phi, t)|^{2}\\
    &\leq \lambda_{3}|\varphi(1)-\phi(1)|^{2}+\lambda_{4}\int^{\overline{\theta}}_{\underline{\theta}} |\varphi(\theta)-\phi(\theta)|^{2}\dif \nu(\theta).
    \end{align*}
\end{enumerate}

\bt
Assume  $(\mathrm{H}3)$ and  $(\mathrm{H}4)$. If  the following conditions hold:
\begin{enumerate}
\item[$1^{\circ}$]
$
\lambda_{1}- 2\lambda_{2}(\lfloor\frac{1}{\underline{\theta}}\rfloor +1)>0.
$
\item[$2^{\circ}$] $\zeta^{*}$ is the unique positive root of the following equation:
 \begin{align}\label{3.25}
\lambda_{1}-\eta= 2\lambda_{2}\bigg(\bigg.\bigg\lfloor\bigg.\frac{1}{\underline{\theta}}\bigg\rfloor\bigg. +1\bigg)\bigg.\underline{\theta}^{-\zeta}.
\end{align}
\item[$3^{\circ}$]$\Delta$  is sufficiently small such that:
\begin{align}
\bigg[\bigg. \lambda_{2}\bigg(\bigg.\bigg\lfloor\bigg.\frac{1}{\underline{\theta}}\bigg\rfloor\bigg. +1\bigg)\bigg.\underline{\theta}^{-\zeta^{*}-1}  +\lambda_{4}\bigg(\bigg.\bigg\lfloor\bigg.\frac{1}{\underline{\theta}}\bigg\rfloor\bigg. +1\bigg)\bigg.\underline{\theta}^{-\zeta^{*}-1}\Delta\bigg]\bigg.\Delta < \frac{1}{2}.
\end{align}
\end{enumerate}
Then for any $\varepsilon\in (0,\frac{\zeta^{*}}{2}),$  there exists a sufficiently small $\Delta^{*} \in (0,1)$ such that the approximate solution $y(k)$ defined by \eqref{neq} satisfies
  \begin{align*}
\limsup_{k\rightarrow \infty }\frac{\log|y(k)|}{\log|(k+1)\Delta|}\leq -\frac{\zeta^{*}}{2}+\varepsilon,a.s.
\end{align*}
and
\begin{align*}
\limsup_{k\rightarrow \infty }\frac{\log\mE[|y(k)|^{2}]}{\log|(k+1)\Delta|}\leq -\zeta^{*}+2\varepsilon,a.s.
\end{align*}
\et

\begin{proof}
By virtue of  \eqref{neq},  we have
\begin{equation}
\begin{split}
|y(k+1)|^{2}&\le |y(k)|^{2}+ |f(y_k,k\Delta)|^{2}\Delta^{2}+2y^{T}(k)f(y_k,k\Delta)\Delta\\
&+|g(y_k,k\Delta)|^{2}\Delta+|g(y_k,k\Delta)|^{2}((\Delta B(k))^{2}-\Delta)
\\
&+2y^{T}(k)g(y_k,k\Delta)\Delta B(k)+2f(y_k,k\Delta)g(y_k,k\Delta)\Delta B(k)\Delta\\
&\leq |y(k)|^{2}+ |f(y_k,k\Delta)|^{2}\Delta^{2}+2y^{T}(k)f(y_k,k\Delta)\Delta\\
&+|g(y_k,k\Delta)|^{2}\Delta+M(k),
\end{split}
\end{equation}
According to $(\mathrm{H}3)$ and $(\mathrm{H}4)$,  we get
\begin{align}\label{3.27}
|y(k+1)|^{2}-|y(k)|^{2}&\leq\bigg(\bigg. -\lambda_{1}|y(k)|^{2}+ \lambda_{2}\int^{\overline{\theta}}_{\underline{\theta}}|y_{k}( k\theta\Delta)|^{2}\dif \nu(\theta)\bigg)\bigg. \Delta \no\\
&+\bigg(\bigg. \lambda_{3}|y(k)|^{2}+ \lambda_{4}\int^{\overline{\theta}}_{\underline{\theta}}|y_{k}( k\theta\Delta)|^{2}\dif \nu(\theta)\bigg)\bigg. \Delta^{2} +M(k).
\end{align}
Multiplying $(1+(1+k)\Delta)^{\gamma}$ on both sides of the inequality \eqref{3.27} yields that
\begin{align}\label{3.28}
&(1+(1+k)\Delta)^{\gamma}|y(k+1)|^{2}-(1+k\Delta)^{\gamma}|y(k)|^{2}\no\\
&\leq(1+(1+k)\Delta)^{\gamma}\bigg(1\bigg.-\frac{(1+k\Delta)^{\gamma})}{(1+(1+k)\Delta)^{\gamma}}\bigg)\bigg.|y(k)|^{2}\no\\
&+\bigg(\bigg. -\lambda_{1}(1+(1+k)\Delta)^{\gamma}|y(k)|^{2}+ \lambda_{2}\int^{\overline{\theta}}_{\underline{\theta}}(1+(1+k)\Delta)^{\gamma}|y_{k}( k\theta\Delta)|^{2}\dif \nu(\theta)\bigg)\bigg. \Delta \no\\
&+\bigg(\bigg. \lambda_{3}(1+(1+k)\Delta)^{\gamma}|y(k)|^{2}+ \lambda_{4}\int^{\overline{\theta}}_{\underline{\theta}}(1+(1+k)\Delta)^{\gamma}|y_{k}( k\theta\Delta)|^{2}\dif \nu(\theta)\bigg)\bigg. \Delta^{2} +M(k),
\end{align}
where $\gamma$ is a positive constant. Observing  that $1-|x|^{\gamma} \leq -\gamma\log|x|,$ one has
$$1-\frac{(1+k\Delta)^{\gamma}}{(1+(1+k)\Delta)^{\gamma}}\leq \gamma\log\frac{1+(1+k)\Delta}{1+k\Delta} \leq \frac{\gamma\Delta}{1+k\Delta}\leq \gamma\Delta.$$
Summing both sides of \eqref{3.28} from $j=0$ to $j=k-1,$  we obtain
\begin{align}\label{3.29}
&(1+k\Delta)^{\gamma}|y(k)|^{2}\no\\
&\leq x_{0}+\sum^{k-1}_{j=0}(1+(j+1)\Delta)^{\gamma}\gamma\Delta|y(j)|^{2}\no\\
&+\bigg(\bigg. -\lambda_{1}\sum^{k-1}_{j=0}(1+(j+1)\Delta)^{\gamma}|y(j)|^{2}+ \lambda_{2}\sum^{k-1}_{j=0}\int^{\overline{\theta}}_{\underline{\theta}}(1+(j+1)\Delta)^{\gamma}|y_{j}( j\theta\Delta)|^{2}\dif \nu(\theta)\bigg)\bigg. \Delta \no\\
&+\bigg(\bigg. \lambda_{3}\sum^{k-1}_{j=0}(1+(j+1)\Delta)^{\gamma}|y(j)|^{2}+ \lambda_{4}\int^{\overline{\theta}}_{\underline{\theta}}\sum^{k-1}_{j=0}(1+(j+1)\Delta)^{\gamma}|y_{j}( j\theta\Delta)|^{2}\dif \nu(\theta)\bigg)\bigg. \Delta^{2}\no\\
&\leq x_{0}+\sum^{k-1}_{j=0}(1+(j+1)\Delta)^{\gamma}\gamma\Delta|y(j)|^{2}\no\\
&+\bigg(\bigg. -\lambda_{1}\sum^{k-1}_{j=0}(1+(j+1)\Delta)^{\gamma}|y(j)|^{2}+ \lambda_{2}\sum^{k-1}_{j=0}\int^{\overline{\theta}}_{\underline{\theta}}(1+(j+1)\Delta)^{\gamma}|y(\lfloor j\theta \rfloor )|^{2}\dif \nu(\theta)\no\\
&\quad \quad \quad\quad \quad \quad\quad \quad \quad\quad \quad \quad + \lambda_{2}\sum^{k-1}_{j=0}\int^{\overline{\theta}}_{\underline{\theta}}(1+(j+1)\Delta)^{\gamma}|y(\lfloor j\theta \rfloor+1)|^{2}\dif \nu(\theta)\bigg)\bigg. \Delta \no\\
&+\bigg(\bigg. \lambda_{3}\sum^{k-1}_{j=0}(1+(j+1)\Delta)^{\gamma}|y(j)|^{2}+ \lambda_{4}\sum^{k-1}_{j=0}\int^{\overline{\theta}}_{\underline{\theta}}(1+(j+1)\Delta)^{\gamma}|y(\lfloor j\theta \rfloor )|^{2}\dif \nu(\theta)\no\\
&\quad \quad \quad\quad \quad \quad\quad \quad \quad\quad \quad \quad + \lambda_{4}\sum^{k-1}_{j=0}\int^{\overline{\theta}}_{\underline{\theta}}(1+(j+1)\Delta)^{\gamma}|y(\lfloor j\theta \rfloor+1)|^{2}\dif \nu(\theta)\bigg)\bigg. \Delta^{2}\no\\
 & +\sum^{k-1}_{j=0}M(j).
\end{align}
Firstly,  we compute
\begin{align}
\sum^{k-1}_{j=0}\int^{\overline{\theta}}_{\underline{\theta}}(1+(j+1)\Delta)^{\gamma}|y( \lfloor j\theta \rfloor )|^{2}\dif \nu(\theta) =\int^{\overline{\theta}}_{\underline{\theta}}\sum^{k-1}_{j=0}(1+(j+1)\Delta)^{\gamma}|y( \lfloor j\theta \rfloor )|^{2}\dif \nu(\theta),
\end{align}
and
\begin{align}
\sum^{k-1}_{j=0}\int^{\overline{\theta}}_{\underline{\theta}}(1+(j+1)\Delta)^{\gamma}|y(\lfloor j\theta \rfloor+1)|^{2}\dif \nu(\theta) =\int^{\overline{\theta}}_{\underline{\theta}}\sum^{k-1}_{j=0}(1+(j+1)\Delta)^{\gamma}|y(\lfloor j\theta \rfloor+1)|^{2}\dif \nu(\theta).
\end{align}
It is obvious that
\begin{align}\label{3.31}
&\sum^{k-1}_{j=0}(1+(j+1)\Delta)^{\gamma}|y(\lfloor j\theta \rfloor )|^{2}\no\\
&=|y(0) |\sum_{0\leq j\leq k-1:\lfloor \theta j\rfloor =0}(1+(j+1)\Delta)^{\gamma}+|y(1) |\sum_{0\leq j \leq k-1:\lfloor \theta j\rfloor =1}(1+(j+1)\Delta)^{\gamma}\no\\
&+\cdots + |y(\lfloor \theta (k-1)\rfloor) |\sum_{0\leq j\leq k-1:\lfloor \theta j\rfloor =\lfloor \theta (k-1)\rfloor}(1+(j+1)\Delta)^{\gamma},
\end{align}
and
\begin{align}\label{3.31+}
\sum^{k-1}_{j=0}&(1+(j+1)\Delta)^{\gamma}|y(\lfloor j\theta \rfloor+1)|^{2}\no\\
&=|y(1) |\sum_{0\leq j \leq k-1:\lfloor \theta j\rfloor =0}(1+(j+1)\Delta)^{\gamma}+|y(2) |\sum_{0\leq j \leq k-1:\lfloor \theta j\rfloor =1}(1+(j+1)\Delta)^{\gamma}\no\\
&+\cdots +|y(\lfloor \theta (k-1)\rfloor+1) |\sum_{0\leq j \leq k-1:\lfloor \theta j\rfloor =\lfloor \theta (k-1)\rfloor}(1+(j+1)\Delta)^{\gamma},
\end{align}
Let $\Delta_{0}=\frac{1-\overline{\theta}}{\overline{\theta}}. $ Since $\overline{\theta}>\frac{1}{2}$,   $\Delta_{0}<1$. Moreover,  and for any $\Delta \in (0, \Delta_{0}), \theta\le \overline{\theta} $, we can derive $\theta+ \theta\Delta\le 1.$ Similar to that of the proof of Theorem \ref{YT1},  \eqref{3.31} and \eqref{3.31+}  can be written as follows, respectively.
\begin{align}
&\sum_{j=0}^{k-1}(1+(j+1)\Delta)^{\gamma}|y(\lfloor \theta j\rfloor)|^{2} \leq \bigg(\bigg.\bigg\lfloor\bigg.\frac{1}{\theta}\bigg\rfloor\bigg. +1\bigg)\bigg.\sum^{\lfloor\theta (k-1)\rfloor}_{j=0}\bigg(\bigg.1+\bigg(\bigg.\frac{j+1}{\theta}+1\bigg)\bigg.\Delta\bigg)\bigg.^{\gamma}|y(j)|^{2}\no\\
&\leq \bigg(\bigg.\bigg\lfloor\bigg.\frac{1}{\theta}\bigg\rfloor\bigg. +1\bigg)\bigg.\theta^{-\gamma}\sum^{\lfloor\theta (k-1)\rfloor}_{j=0}(\theta(1+\Delta)+(j+1)\Delta)^{\gamma}|y(j)|^{2}\no\\
&\leq \bigg(\bigg.\bigg\lfloor\bigg.\frac{1}{\theta}\bigg\rfloor\bigg. +1\bigg)\bigg.\theta^{-\gamma}\sum^{k-1}_{j=0}(1+(j+1)\Delta)^{\gamma}|y(j)|^{2}\no\\
&\leq \bigg(\bigg.\bigg\lfloor\bigg.\frac{1}{\underline{\theta}}\bigg\rfloor\bigg. +1\bigg)\bigg.\underline{\theta}^{-\gamma}\sum^{k-1}_{j=0}(1+(j+1)\Delta)^{\gamma}|y(j)|^{2},
\end{align}
and
\begin{align}
&\sum_{j=0}^{k-1}(1+(j+1)\Delta)^{\gamma}|y(\lfloor \theta j\rfloor+1)|^{2} \leq \bigg(\bigg.\bigg\lfloor\bigg.\frac{1}{\theta}\bigg\rfloor\bigg. +1\bigg)\bigg.\sum^{\lfloor\theta (k-1)\rfloor}_{j=0}\bigg(\bigg.1+\bigg(\bigg.\frac{j+1}{\theta}+1\bigg)\bigg.\Delta\bigg)\bigg.^{\gamma}|y(j+1)|^{2}\no\\
&\leq \bigg(\bigg.\bigg\lfloor\bigg.\frac{1}{\underline{\theta}}\bigg\rfloor\bigg. +1\bigg)\bigg.\underline{\theta}^{-\gamma}\sum^{k}_{j=1}(1+j\Delta)^{\gamma}|y(j)|^{2}.
\end{align}
Combining with $\eqref{3.29},$  we arrive at
 \begin{align}
&(1+k\Delta)^{\gamma}|y(k)|^{2}\no\\
&\leq x_{0}+\sum^{k-1}_{j=0}(1+(j+1)\Delta)^{\gamma}\gamma\Delta|y(j)|^{2}\no\\
&+\bigg[\bigg. -\lambda_{1}\sum^{k-1}_{j=0}(1+(j+1)\Delta)^{\gamma}|y(j)|^{2}+ \lambda_{2}\bigg(\bigg.\bigg\lfloor\bigg.\frac{1}{\underline{\theta}}\bigg\rfloor\bigg. +1\bigg)\bigg.\underline{\theta}^{-\gamma}\sum^{k-1}_{j=0}(1+(j+1)\Delta)^{\gamma}|y(j)|^{2}\no\\
&\quad \quad\quad\quad \quad\quad\quad \quad\quad\quad \quad\quad+ \lambda_{2}\bigg(\bigg.\bigg\lfloor\bigg.\frac{1}{\underline{\theta}}\bigg\rfloor\bigg. +1\bigg)\bigg.\underline{\theta}^{-\gamma}\sum^{k}_{j=1}(1+j\Delta)^{\gamma}|y(j)|^{2} \bigg]\bigg. \Delta \no\\
&+\bigg[\bigg. \lambda_{3}\sum^{k-1}_{j=0}(1+(j+1)\Delta)^{\gamma}|y(j)|^{2}+ \lambda_{4}\bigg(\bigg.\bigg\lfloor\bigg.\frac{1}{\underline{\theta}}\bigg\rfloor\bigg. +1\bigg)\bigg. \underline{\theta}^{-\gamma}\sum^{k-1}_{j=0}(1+(j+1)\Delta)^{\gamma}|y(j)|^{2}\no\\
&\quad \quad\quad\quad \quad\quad\quad \quad\quad\quad \quad\quad+ \lambda_{4}\bigg(\bigg.\bigg\lfloor\bigg.\frac{1}{\underline{\theta}}\bigg\rfloor\bigg. +1\bigg)\bigg.\underline{\theta}^{-\gamma}\sum^{k}_{j=1}(1+j\Delta)^{\gamma}|y(j)|^{2} \bigg]\bigg. \Delta^{2}\no\\
& +\sum^{k-1}_{j=0}M(j)\no\\
&\leq x_{0}+\bigg[\bigg. -\lambda_{1}+\gamma+ \lambda_{2}\bigg(\bigg.\bigg\lfloor\bigg.\frac{1}{\underline{\theta}}\bigg\rfloor\bigg. +1\bigg)\bigg.\underline{\theta}^{-\gamma}
 +\lambda_{3}\Delta\no\\
  &\quad \quad\quad\quad \quad\quad\quad \quad\quad\quad \quad\quad+\lambda_{4}\bigg(\bigg.\bigg\lfloor\bigg.\frac{1}{\underline{\theta}}\bigg\rfloor\bigg. +1\bigg)\bigg.\underline{\theta}^{-\gamma}\Delta\bigg]\bigg.\sum^{k-1}_{j=0}(1+(j+1)\Delta)^{\gamma}|y(j)|^{2} \Delta\no\\
&+\bigg[\bigg. \lambda_{2}\bigg(\bigg.\bigg\lfloor\bigg.\frac{1}{\underline{\theta}}\bigg\rfloor\bigg. +1\bigg)\bigg.\underline{\theta}^{-\gamma}  +\lambda_{4}\bigg(\bigg.\bigg\lfloor\bigg.\frac{1}{\underline{\theta}}\bigg\rfloor\bigg. +1\bigg)\bigg.\underline{\theta}^{-\gamma}\Delta\bigg]\bigg.\sum^{k}_{j=1}(1+j\Delta)^{\gamma}|y(j)|^{2} \Delta\no\\
&+\sum^{k-1}_{j=0}M(j)\no\\
&\leq x_{0}+\bigg[\bigg. -\lambda_{1}+\gamma+ 2\lambda_{2}\bigg(\bigg.\bigg\lfloor\bigg.\frac{1}{\underline{\theta}}\bigg\rfloor\bigg. +1\bigg)\bigg.\underline{\theta}^{-\gamma}
 +\lambda_{3}\Delta\no\\
 &\quad \quad\quad\quad \quad\quad\quad \quad\quad\quad \quad\quad+2\lambda_{4}\bigg(\bigg.\bigg\lfloor\bigg.\frac{1}{\underline{\theta}}\bigg\rfloor\bigg. +1\bigg)\bigg.\underline{\theta}^{-\gamma}\Delta\bigg]\bigg.\sum^{k-1}_{j=0}(1+(j+1)\Delta)^{\gamma}|y(j)|^{2} \Delta\no\\
&+\bigg[\bigg. \lambda_{2}\bigg(\bigg.\bigg\lfloor\bigg.\frac{1}{\underline{\theta}}\bigg\rfloor\bigg. +1\bigg)\bigg.\underline{\theta}^{-\gamma}  +\lambda_{4}\bigg(\bigg.\bigg\lfloor\bigg.\frac{1}{\underline{\theta}}\bigg\rfloor\bigg. +1\bigg)\bigg.\underline{\theta}^{-\gamma}\Delta\bigg]\bigg.(1+k\Delta)^{\gamma}|y(k)|^{2} \Delta\no\\
&+\sum^{k-1}_{j=0}M(j).
\end{align}
 Set
 \begin{align}
H(\gamma, \Delta)= -\lambda_{1}+\gamma+ 2\lambda_{2}\bigg(\bigg.\bigg\lfloor\bigg.\frac{1}{\underline{\theta}}\bigg\rfloor\bigg. +1\bigg)\bigg.\underline{\theta}^{-\gamma}
 +\lambda_{3}\Delta
 +2\lambda_{4}\bigg(\bigg.\bigg\lfloor\bigg.\frac{1}{\underline{\theta}}\bigg\rfloor\bigg. +1\bigg)\bigg.\underline{\theta}^{-\gamma}\Delta.
\end{align}
Immediately, one can see that
\begin{align*}
\frac{\dif H(\gamma, \Delta)}{\dif \gamma}= 1+ 2\lambda_{2}\bigg(\bigg.\bigg\lfloor\bigg.\frac{1}{\underline{\theta}}\bigg\rfloor\bigg. +1\bigg)\bigg.\underline{\theta}^{-\gamma}\log\frac{1}{\underline{\theta}}+2\lambda_{4}\bigg(\bigg.\bigg\lfloor\bigg.\frac{1}{\underline{\theta}}\bigg\rfloor\bigg. +1\bigg)\bigg.\underline{\theta}^{-\gamma}\log\frac{1}{\underline{\theta}}\Delta>0,
\end{align*}
and
\begin{align*}
H(0, \Delta)=-\lambda_{1}+ 2\lambda_{2}\bigg(\bigg.\bigg\lfloor\bigg.\frac{1}{\underline{\theta}}\bigg\rfloor\bigg. +1\bigg)\bigg.
+\lambda_{3}\Delta +2\lambda_{4}\bigg(\bigg.\bigg\lfloor\bigg.\frac{1}{\underline{\theta}}\bigg\rfloor\bigg. +1\bigg)\bigg.\Delta<0,
\end{align*}
for $\Delta \in (0, \Delta_{1}\wedge \Delta_{0}), \Delta_{1}=\frac{\lambda_{1}- \lambda_{2}(\lfloor\frac{1}{\underline{\theta}}\rfloor +1)}{\lambda_{3}+\lambda_{4}(\lfloor\frac{1}{\underline{\theta}}\rfloor +1)}.$
Then, for any $\Delta \in  (0,  \Delta_{1}\wedge \Delta_{0})), $ it follows that there exists a constant $\gamma^{*}_{\Delta}$ such that $H(\gamma^{*}_{\Delta}, \Delta)=0.$
This together with condition $3^{\circ}$ implies
\begin{align}\label{3.35}
(1+k\Delta)^{\gamma^{*}_{\Delta}}|y(k)|^{2}
\leq 2|x_{0}|^{2}+ 2\sum^{k-1}_{j=0}M(j).
\end{align}
Since, $\sum^{k-1}_{j=0}M(j)$ is a martingale, we get
\begin{align}\label{3.36+}
\limsup_{k\rightarrow \infty }(1+k\Delta)^{\gamma^{*}_{\Delta}}\mE|y(k)|^{2}
< \infty.
\end{align}
From Lemma 2.1, we get
\begin{align}\label{3.36}
\limsup_{k\rightarrow \infty }(1+k\Delta)^{\gamma^{*}_{\Delta}}|y(k)|^{2}
< \infty.
\end{align}
Noting that
 $$ \lim_{\Delta \rightarrow 0}H(\gamma, \Delta)=-\lambda_{1}+\gamma+ 2\lambda_{2}\bigg(\bigg.\bigg\lfloor\bigg.\frac{1}{\underline{\theta}}\bigg\rfloor\bigg. +1\bigg)\bigg.\underline{\theta}^{-\gamma},$$
and \eqref{3.25},  one has $\lim_{\Delta \rightarrow 0}\gamma^{*}_{\Delta}= \zeta^{*}.$
Thus, for any $\varepsilon \in  (0, \frac{\zeta^{*}}{2}),$ there exists $\Delta_{2}$ such that $\gamma^{*}_{\Delta}>\zeta^{*} -2\varepsilon$ for any $\Delta \in (0,\Delta_{2}).$
Then, for any $\Delta \in (0, \Delta_{0}\wedge \Delta_{1}\wedge \Delta_{2}),$  \eqref{3.35}, \eqref{3.36+} and  \eqref{3.36} imply the result in the theorem.
The proof is complete.
\end{proof}
\br
Under the conditions in Theorem 4.3,  from Theorem 3.6 in \cite{why}, we know that the analytical solution of \eqref{eq41} has the property of polynomial stability. This shows that
the EM method inherits the polynomial stability of the true solution.
\er
\begin{exa}
{\rm
Consider the following equation:
\begin{align}\label{EY2}
 &\dif x(t)= f( x_t, t)\dif t +g(  x_t, t)\dif  B(t), t\in [0, \infty)
 & x(0)=x_{0},
 \end{align}
where
\begin{align*}
f( \varphi, t)=
-0.4 \varphi(1)+0.04\int^{\frac{4}{5}}_{\frac{3}{4}}|\varphi(\theta )|\dif \nu (\theta).
\end{align*}
and
\begin{align*}
g( \varphi, t)=
0.3\int^{\frac{4}{5}}_{\frac{3}{4}}|\varphi(\theta )|\dif \nu (\theta).
\end{align*}
Form the above definition,   it follows that
\begin{align*}
&2\langle \varphi(1)-\phi(1), f(\varphi, t)-f(\phi,t)\rangle +|g(\varphi,t)-g(\phi,t)|^{2}\\
&=2\langle \varphi(1)-\phi(1), -0.4(\varphi(1)-\phi(1))+0.04\int^{\frac{4}{5}}_{\frac{3}{4}}(\varphi(\theta )-\phi(\theta ))\dif \nu (\theta)\rangle \\
&+\bigg|\bigg.0.3\int^{\frac{4}{5}}_{\frac{3}{4}}(\varphi(\theta t)-\phi(\theta t))\dif \nu (\theta)\bigg|\bigg.^{2}\\
& \leq -0.8|\varphi(1)-\phi(1)|^{2}+0.08(\varphi(1)-\phi(1))\int^{\frac{4}{5}}_{\frac{3}{4}}(\varphi(\theta )-\phi(\theta ))\dif \nu (\theta)\\
&+\bigg|\bigg.0.3\int^{\frac{4}{5}}_{\frac{3}{4}}(\varphi(\theta )-\phi(\theta ))\dif \nu (\theta)\bigg|\bigg.^{2}\\
& \leq -0.76|\varphi(1)-\phi(1)|^{2}+0.13\int^{\frac{4}{5}}_{\frac{3}{4}}|(\varphi(\theta )-\phi(\theta ))|^{2}\dif \nu (\theta),
\end{align*}
and
\begin{align*}
& |f(\varphi, t)-f(\phi,t)|^{2}\vee |g(\varphi, t)-g(\phi,t)|^{2}\\
& \leq 0.19|\varphi(1)-\phi(1)|^{2}+0.09\int^{\frac{4}{5}}_{\frac{3}{4}}|(\varphi(\theta )-\phi(\theta ))|^{2}\dif \nu (\theta).
\end{align*}
Letting $\lambda_{1}= 0.76,\lambda_{2}= 0.13, \lambda_{3}= 0.19,\lambda_{1}= 0.09, \underline{\theta}=\frac{3}{4},$
it can be  seen that the conditions in Theorem 3.6 are satisfied.
We therefore conclude that the numerical solutions of \eqref{EY2} are almost surely polynomial stable. }

\end{exa}

\end{document}